\documentclass[reqno]{amsart} 
\usepackage{amsfonts, amsmath, amsthm, amssymb, latexsym, graphicx, geometry, xcolor, colortbl,  mathtools}

\theoremstyle{plain}
\newtheorem{theorem}{Theorem}[section]
\newtheorem{corollary}[theorem]{Corollary}
\newtheorem{lemma}[theorem]{Lemma}
\newtheorem{proposition}[theorem]{Proposition}
\theoremstyle{definition}
\newtheorem{definition}[theorem]{Definition}
\newtheorem{remarks}[theorem]{Remarks}
\newtheorem{remark}[theorem]{Remark}

\numberwithin{equation}{section}

\title[free multiplicative Brownian motions]{An explicit formula for free multiplicative Brownian motions via spherical functions}

\author{Martin Auer, Michael Voit} 
\address{Fakult\"at Mathematik, Technische Universit\"at Dortmund,
          Vogelpothsweg 87,
          D-44221 Dortmund, Germany}
\email{martin.auer@math.tu-dortmund.de, michael.voit@math.tu-dortmund.de}


\begin{document}
\date{\today}

\begin{abstract} After some normalization, the logarithms of the ordered singular values of Brownian motions on $GL(N,\mathbb F)$ 
  with   $\mathbb F=\mathbb R, \mathbb C$ form Weyl-group invariant Heckman-Opdam processes on $\mathbb R^N$ of type $A_{N-1}$.
  We use classical elementary formulas for the
  spherical functions of   $GL(N,\mathbb C)/SU(N)$ and the associated  Euclidean spaces  $H(N,\mathbb C)$ of Hermitian matrices,
  and show that  in the   $GL(N,\mathbb C)$-case, these processes can be
also interpreted as 
ordered eigenvalues of Brownian motions on $H(N,\mathbb C)$ with particular drifts.  This leads to an explicit
description for the free limits for the associated empirical processes for $N\to\infty$ where
these limits are independent from the parameter $k$ of the Heckman-Opdam processes. In particular we get  new formulas for the distributions
of the free multiplicative Browniam motion of Biane.
We also show how this approach works for the root systems $B_N, C_N, D_N$.
\end{abstract}

\maketitle

\centerline{ This paper is dedicated to Tom Koornwinder at the occasion of his 80th birthday}

\section{Introduction}
For  $\mathbb F=\mathbb R, \mathbb C$, and the quaternions $ \mathbb H$, consider the
right-Brownian motions $(G_t)_{t\ge0}$ on the 
Lie groups $GL(N,\mathbb F)$ as solutions of the stochastic differential equations $\partial G_t=(\partial W_t)\cdot G_t$
in the Stratonovitch sense with start in the identity matrix $G_0=I_N$  where $(W_t)_{t\geq0}$ is a Brownian motion on $\mathbb{F}^{N\times N}$.
It is well known (see e.g. \cite{NRW})
that  the processes $(G_t\bar G_t^T)_{t\ge0}$ then are diffusions
on the cones $P(N,\mathbb F)$  of positive definite matrices such that the processes of their ordered  eigenvalues
$(\lambda_t=(\lambda_{t,1},\ldots,\lambda_{t,N}))_{t\geq0}$ on the Weyl chamber
$\{\lambda\in ]0,\infty[^N: \>  \lambda_{1}\le \ldots\le \lambda_{N}\}$ never collide for $t>0$ almost surely. Furthermore, their
(normalized)  logarithms
$$\left(X_t:=\frac{1}{2}\ln \lambda_t:=\left(\frac{1}{2}\ln \lambda_{t,1},\ldots,\frac{1}{2}\ln \lambda_{t,N}\right)\right)_{t\ge0}$$
  satisfy the Ito-type SDEs
  \begin{equation}\label{SDE1}
    dX_{t,i}= dB_{t,i}+ \frac{d}{2}\sum_{j: j\ne i} \coth(X_{t,i}-X_{t,j}) dt \quad\quad(i=1,\ldots, N)\end{equation}
  on the closed Weyl chambers $C_N^A:=\{x\in\mathbb R^N: \>  x_1\le\ldots\le x_N\}$
 of type A with $d:=dim_{\mathbb R}\mathbb F\in\{1,2,4\}$ and initial condition $X_{t,0}=0\in C_N^A$
  where the  $(B_{t,i})_{t\ge0}$ are independent Brownian motions, and  the $(X_t)_{t\ge0}$ are in the interior of $C_N^A$ a.s.
  for $t>0$.

  The SDEs (\ref{SDE1}) can be embedded into the general framework of the SDEs
   \begin{equation}\label{SDE2}
    dX_{t,i}^{(k)}= dB_{t,i}+ k\sum_{j\colon j\ne i} \coth\left(X^{(k)}_{t,i}-X^{(k)}_{t,j}\right) dt \quad\quad(i=1,\ldots, N)\end{equation}
   on  $C_N^A$ with $X^{(k)}_{t,0}=0\in C_N^A$  for arbitrary ``coupling constants'' $k\in ]0,\infty[$
    which appear in the context of integrable particle models
    of Calogero-Moser-Sutherland type; see e.g. \cite{DV, F} for the background.
 The SDEs (\ref{SDE2}) have unique strong solutions for all $k\ge 1/2$
    in the interior of $C_N^A$ for $t>0$ a.s. for all starting points in  $C_N^A$;
    see e.g. \cite{GM} for general results in this direction.
    The solutions  $(X_{t}^{(k)})_{t\geq0}$ are also known as Weyl-group invariant Heckman-Opdam processes
    related to the root systems of type A from a more analytic point of view; see \cite{Sch1, Sch2}.
    Their transition probabilities can be expressed in terms of inverse spherical Fourier-type transforms related to
    Heckman-Opdam hypergeometric functions of  type A; for the background on these functions we refer to the monographs \cite{HS, HO}.
    
For our purpose we also consider the time-transformed processes
$\left(\tilde{X}_t^{(k),N}:=\tilde X_{t}^{(k)}:=X_{t/k}^{(k)}\right)_{t\ge0}$ which satisfy
\begin{equation}\label{SDE3}
  d\tilde{X}_{t,i}^{(k)}= \frac{1}{\sqrt k}dB_{t,i}+ \sum_{j: j\ne i} \coth(\tilde{X}_{t,i}^{(k)}-\tilde{X}_{t,j}^{(k)}) dt \quad\quad(i=1,\ldots, N).\end{equation}
For $k=\infty$, these SDEs formally degenerate  into the ODEs
\begin{equation}\label{ODE}
  \frac{d}{dt} \tilde{X}_{t,i}^{(\infty)}=\sum_{j: j\ne i} \coth(\tilde{X}_{t,i}^{(\infty)}-\tilde{X}_{t,j}^{(\infty)})  \quad\quad(i=1,\ldots, N).\end{equation}
By the methods in \cite{VW2, AVW} for similar situations, the ODEs (\ref{ODE})  admit unique solutions for
all initial values in $C_N^A$ where these solution  are in the interior of  $C_N^A$ for $t>0$.

In \cite{B1, B2}, Biane  introduced the free multiplicative Brownian motion $(g_t)_{t\geq0}$ as  solution of the free SDE $dg_t=idc_tg_t$ with
 a circular Brownian motion $c_t$. 
 For $\mathbb{F}=\mathbb{C}$ it was shown in \cite{Ke1} that in the large $N$ limit the
 right-Brownian motions  $(G_t)_{t\geq0}$ converge as  noncommutative stochastic processes to $(g_t)_{t\geq0}$.
In particular the spectral distribution of $G_t\overline{G_t}^T$ converges a.s. to the spectral distribution $\mu_t$ of $g_tg_t^*$ where
$(\mu_t)_{t\geq0}$ forms a free multiplicative convolution semigroup, and for each $t$ the measure $\mu_t$
is compactly supported on $(0,\infty)$ and absolutely continuous w.r.t. Lebesgue measure.
However, the density of $\mu_t$ can only be stated implicitly, see section 4.2.5 in \cite{B2}.
In the literature $g_tg_t^*$ is often referred to as 'free positive multiplicative Brownian motion'
since it has free positive multiplicative increments; see \cite{B3} for the latter terminology.
Further properties of its spectral measure can e.g. be found in \cite{Z,HUW}.

Note that the spectral distribution of $G_tG_t^*$ is just the image measure of the empirical distribution of $X_t^{(1)}$ under the map $x\mapsto e^{2x}$.
It is expected that the convergence of the empirical measures for $N\to\infty$ extends to general $k\in[1/2,\infty]$.
For $k=1/2$, corresponding to $\mathbb{F}=\mathbb{R}$, this was shown in Section 7 of \cite{MP}, where $k$ just enters as a time scaling.
In this paper we are going to show the extension of this convergence to general $k\in[1/2,\infty]$, and, as a main result,
the following identification of the limit measure via the free additive convolution of a uniform distribution $U_r$ on $[-r,r]$
with a semicircle distribution $\mu_{sc,s}\in M^1(\mathbb R)$ with some radius $s>0$ which has the Lebesgue density $$\frac{2}{\pi s^2}\sqrt{s^2-x^2}{\bf 1}_{[-s,s]}(x).$$

\begin{theorem}\label{main-theorem1}
  For each $N\in \mathbb N$ and  $k\in [1/2,\infty]$ consider the solutions $(\tilde X_{t}^{(k)})_{t\geq0}$ of (\ref{SDE3}) and
  (\ref{ODE}) for $k<\infty$ and $k=\infty$ respectively with start in $0\in C_N^A$.
Then, for all $k$ and $t>0$, the  normalized empirical measures
	\begin{equation}
		\mu_{N,t}:=\frac{1}{N}\sum_{i=1}^N\delta_{\tilde{X}_{t/N,i}^{(k)}}
	\end{equation}
 tend weakly to $U_{t}\boxplus\mu_{sc,2\sqrt{t}}$ almost surely.
\end{theorem}

Further details and the supports of these limit measures can be found  in  \cite{B2, HZ, Ke2}.
The fact that the limit measure behaves for small $t$ like a semicircle distribution and for large $t$
like a uniform distribution can also be seen by the corresponding asymptotics for $\tilde X_{t}^{(1)}$  in \cite{IS}.
However, our precise description of the transition between these two regimes via free additive convolutions seems to be new.

Our proof of Theorem \ref{main-theorem1} decomposes into two parts. In the first part one shows that all moments of the empirical measures 
$\mu_{N,t}$ converge weakly  almost surely to the moments of some probability distribution which does not depend
on $k\in [1/2,\infty]$.
By the moment convergence theorem, this implies weak convergence of the measures a.s. This part is similar to other models
like for Dyson Brownian motions as e.g.~in Ch.~4 of \cite{AGZ} and multivariate Bessel and Jacobi processes in \cite{VW1, AVW}.
For the identification of the  limits in a second step,
we   study the case $k=1$ which corresponds to the groups $GL(N,\mathbb C)$.
In this complex case there is a  classical connection between the spherical functions of the Gelfand pairs
$(GL(N,\mathbb C), SU(N))$ and  $(SU(N)\ltimes H_N(\mathbb C), SU(N))$, where $SU(N)\ltimes H_N(\mathbb C)$ is the semidirect product when 
$SU(N)$ acts on the vector space $H_N(\mathbb C)$
of all $N$-dimensional complex Hermitian matrices as usual by the conjugation $(U,H)\mapsto \bar U H U$; see  Helgason \cite{He}. 
This  connection was used in \cite{Kl2, RV2} to derive some algebraic connection between $SU(N)$-biinvariant
random walks  on  $GL(N,\mathbb C)$
and  special random walks on  $H_N(\mathbb C)$ with  drifts. This connection was also used in \cite{RV2} to derive  limit theorems for
$SU(N)$-biinvariant  random walks  on  $GL(N,\mathbb C)$ for time $t\to\infty$ in a simple way from
corresponding classical limit theorems on the vector spaces
$H_N(\mathbb C)$. We shall apply this connection  to
Brownian motions on  $GL(N,\mathbb C)$  and Brownian motions on $H_N(\mathbb C)$ with particular drifts.
To describe this connection,
we first consider some Brownian motion $(B_t^H)_{t\ge0}$ on $H_N(\mathbb C)$ with one-dimensional Brownian motions on the diagonal entries and
$(j,l)$-entries of the form $\frac{1}{\sqrt 2}(B_{t,j,l}^{\operatorname{Re}}+iB_{t,j,l}^{\operatorname{Im}})$ for $j<l$ with one-dimensional  Brownian motions
$B_{j,l}^{\operatorname{Re}}, B_{j,l}^{\operatorname{Im}}$ where all processes are independent with start in $0$. It is well known that  for $t>0$,
the ordered eigenvalues in $C_N^A$ here have the densities
\begin{equation}\label{density-A-1}
\frac{1}{(2\pi)^{N/2} 2!3!\cdots (N-1)! t^{N^2/2}} e^{-\|x\|^2/(2t)} \prod_{i<j}(x_i-x_j)^{2}\,,
\end{equation}
where $||x\|$ is the usual Euclidean norm on $\mathbb C^N$; see e.g. \cite{AGZ}. Now consider the particular vector
\begin{equation}\label{rho-a}
  \rho:=(-N+1,-N+3,\ldots, N-3,N-1)\in \mathbb R^N
    \end{equation}
which is known in harmonic analysis as weighted half-sum of positive roots (see \cite{BO, He, HS, HO}). Moreover, let $c>0$ be some constant.
 We shall show the following result for the
distributions of the  ordered eigenvalues  of the Brownian motions  $(B_t^H+ tc\cdot \operatorname{diag}(\rho))_{t\ge0}$ on $H_N(\mathbb C)$ with drift
$c\cdot \operatorname{diag}(\rho)\in  H_N(\mathbb C)$.

\begin{proposition}\label{density-a-general-drift}
For $t>0$,
the ordered eigenvalues of  $(B_t^H+ tc\cdot \operatorname{diag}(\rho))_{t\ge0}$ in $C_N^A$ have the density
\begin{equation}\label{density-A-2}
f_{N,c,t}(x):=\frac{e^{-c\|\rho\|^2 t/2}}{c^{N(N-1)/2} (2\pi)^{N/2} 2!3!\cdots (N-1)! t^{N^2/2}} e^{-\|x\|^2/(2t)} \prod_{i<j}((x_i-x_j)\sinh(c(x_i-x_j)))\,,
\end{equation}
with $\lvert\lvert\rho\rvert\rvert^2=(N-1)N(N+1)/3$. Here $\lvert\lvert\cdot\rvert\rvert$ denotes the usual Euclidean norm in $\mathbb{R}^N$.
\end{proposition}

\begin{remark}
	In \cite{TK} the related process $(W_t+(a+t\sigma)\rho)_{t\geq0}$ is considered, where $(W_t)_{t\geq0}$ is a standard real $N$-dimensional Brownian motion and $a,\sigma>0$ are constants.
	Conditioned on the event that the process never hits the boundary of $C_N^{A}$, the corresponding finite dimensional distributions are obtained in \cite{TK} Proposition 1.
	By specializing the latter result to the single-time marginal distribution and taking the limit $a\downarrow0$, one recovers \eqref{density-A-2}.\\
	The fact that the density $f_{N,c,t}$ solves the Fokker-Planck equation corresponding to the SDE \eqref{SDE1} is known by \cite{IS} section 3.
\end{remark}

On the other hand, the following holds for $GL(N,\mathbb C)$:

\begin{proposition}\label{density-a-special-drift}
  Consider the solution  $(X_t)_{t\ge0}$ of the SDE (\ref{SDE1}) with $d=2$ and
  start in $0$ which is related to a Brownian motion on $GL(N,\mathbb C)$ as described above. Then, with the notation from
  Proposition \ref{density-a-general-drift}, $X_t$ has the density $f_{N,1,t}$ for $t>0$.
\end{proposition}

A combination of these two propositions easily leads to Theorem \ref{main-theorem1} where the simple form of the vectors $\rho\in\mathbb R^N$
is responsible for the uniform distributions  
in the  free limit in Theorem \ref{main-theorem1}.

This  connection between Brownian motions on
$GL(N,\mathbb C)$ and $H_N(\mathbb C)$ for the root systems $A_{N-1}$  can be also used for Brownian motions
on  the further  series of 
complex, noncompact  connected simple  Lie groups with finite center which are listed e.g.~in Appendix C in  \cite{Kn}.
We shall discuss these 
results for the groups $SO(2N+1,\mathbb C)$, $G=Sp(N,\mathbb C)$, and $SO(2N,\mathbb C)$
 and the root systems $B_N, C_N, D_N$ respectively briefly in Section 4.

We finally point out that our approach does not seem to work for the  root systems of types $BC_N$ with general parameters
as here in general no such close connection  exists between
the spherical functions on the corresponding symmetric spaces and their Euclidean counterparts.
In fact, in these cases, the  candidates for symmetric spaces of interest are the noncompact Grassmann manifolds over $\mathbb C$.
In these cases there
exist  determinantal formulas for the spherical functions due to Berezin and  Karpelevich \cite{BK}
(see also \cite{Ho} and references there)
in terms of one-dimensional Jacobi functions which extend the spherical functions of the rank-one case; see  \cite{Ko}
for  details on these functions. Moreover, there are related determinantal formulas for the spherical functions
in the associated Euclidean cases where there one-dimensional Bessel functions instead of Jacobi functions appear.
As the connections between  Bessel and  Jacobi functions are not as simple as in the cases mentioned above even for $N=1$
(except for a few trivial cases),
we do not know whether it is possible to derive results for noncompact Grassmann manifolds over $\mathbb C$,
which are similar to those in this paper.

The paper is organized as follows: In Section 2 we  prove Propositions \ref{density-a-general-drift} and \ref{density-a-special-drift}
where these results are embedded into some more general context of projections of Brownian motions with drifts and spherical functions of
complex, noncompact connected semisimple Lie groups.
These results and a corresponding general a.s.~free limit result for the solutions of (\ref{SDE2}),
which involves free multiplicative Brownian motions of Biane, then will lead to 
 Theorem \ref{main-theorem1} in Section 3.
 Finally, in Section 4 we use the ideas of Section 2, in order to derive 
 corresponding results
 for the root systems of types B, C, and D. This leads to connections between new
 free limits of noncompact  multivariate Jacobi processes and associated multivariate Bessel processes studied in
 \cite{ RV1, CGY, AVW, VW1}.

 In general, for SDEs and free probability we refer e.g.~to the textbooks \cite{P} and \cite{NS} respectively.

\section{Spherical functions and Brownian motions}

In this section we derive Propositions \ref{density-a-general-drift} and \ref{density-a-special-drift}. We
start with a quite general setting for general  Euclidean spaces, and only later
 on we  restrict our attention to the setting of Propositions \ref{density-a-general-drift} and \ref{density-a-special-drift}.
The ideas in the following subsection are taken from \cite{V, Kl1, Kl2, RV2}.

\subsection{Projections of invariant processes on Euclidean spaces  with drifts}\label{sub-21}
Let 
$(V, \langle\,.\,,\,.\,\rangle)$ be a 
Euclidean vector space of finite dimension $n$ and
$K\subset O(V)$ a  compact subgroup  of the orthogonal group of $V$.
 Consider the exponential function
 $$e_\lambda(x):=e^{\langle \lambda,x\rangle}  \quad\quad(x\in V).$$
   Furthermore, consider the  space
$$M_e(V):=\{\mu\in M_b(V):\> e_\lambda\cdot \mu\in M_b(V) \quad\text{for all}\quad \lambda\in V\}$$
   of all bounded signed measures on $V$ with exponential moments as well as the subspace
$$M_{e,K}(V):=\{\mu\in M_e(V):\> \mu \quad K\text{-invariant}\}$$
   of those measures which are invariant under the action of $K$ on $V$ in addition.
As the  $e_\lambda$ are multiplicative on $V$, the vector spaces $M_e(V)$ and $M_{e,K}(V)$ are  subalgebras of 
the Banach-algebra $M_b(V)$ with the usual convolution as product where  $M_b(V)$ carries the total variation norm.
Moreover, as
\begin{equation}\label{convo-equation}
  e_\lambda\mu * e_\lambda\nu =  e_\lambda(\mu*\nu) \quad\text{for all}\quad \mu,\nu\in M_e(V), \> \lambda\in V,
  \end{equation}
also the spaces
$$M_{e,K}^\lambda(V):=\{e_\lambda\mu:\> \mu\in M_{e,K}(V) \}$$
are subalgebras of  $M_e(V)$ for all $\lambda\in V$.

Let $V^K:=\{K.x: \> x\in V\}$ be the space of all $K$-orbits in $V$ which is a locally compact Hausdorff space w.r.t.~the quotient topology.
Let $p:V\to V^K$ be the canonical projection, which is continuous and open, 
and  denote the associated push forward of measures also by $p$. We then have the following simple observation; see e.g.~\cite{RV2}:

\begin{lemma}\label{modified-projection} Let  $\lambda\in V$. Then, for $x\in V$, and the normalized Haar measure $dk$ of $K$,
  $$\alpha_\lambda(p(x)):=\int_K e_\lambda(k(x)) \> dk$$
  is a well-defined continuous function on $V^K$, and for all $\mu\in M_{e,K}(V)$,
  \begin{equation}\label{general-image}
    p(e_\lambda \cdot \mu) = \alpha_\lambda \cdot p(\mu).\end{equation}
\end{lemma}

\begin{proof} We only have to check (\ref{general-image}). For this we observe that, by the invariance of $\mu$,
  for each bounded continuous function $f\in C_b (V^K)$,
  \begin{align}
    \int_ {V^K}f \> dp((e_\lambda \mu)&= \int_ {V}f(p(x))e_\lambda( x) \> d\mu(x)=  \int_ {K}\int_ {V} f(p(k(x))e_\lambda( k(x)  )\> d\mu(x)\> dk \notag\\
    &= \int_ {V} f(p(x)) \int_ {K}e_\lambda( k(x)  )\>\> dk \>  d\mu(x) = \int_ {V} f(p(x))\alpha_\lambda(p(x))\> d\mu(x) \notag\\
    &=\int_{V^K} f\cdot  \alpha_ \lambda  \> dp ( \mu) . \notag
\end{align} 
\end{proof}

We next transfer the convolution $*$ on the algebras  $M_{e,K}(V)$ and $M_{e,K}^\lambda(V)$ to convolutions of
 the push forwards of these measures on $V^K$.
We here  follow \cite{RV2}, where measures with compact supports instead of measures with exponential moments are used.
We equip the vector space
$$ M_e(V^K):= \{p(\mu): \mu\in M_{e,K}(V)\}$$
with the total variation norm. Then also
$$ M_e(V^K)= \{p(\mu): \mu\in M_{e,K}^\lambda(V)\},$$
and the mappings
$$ M_{e,K}(V)\to M_e(V^K), \> \mu\mapsto p(\mu) \quad \text{and}\quad M_{e,K}^\lambda(V)\to M_e(V^K), \> \mu\mapsto p(\mu)$$
are isometric isomorphisms of normed spaces. We now transfer the convolution $*$ on the algebras  $M_{e,K}(V)$ and $M_{e,K}^\lambda(V)$ respectively to
$M_e(V^K)$ such that $\mu\mapsto p(\mu)$ in both cases becomes an  isometric isomorphism of normed algebras. Denote these new convolutions
by $\bullet$ and $\bullet_\lambda$ respectively. It can be easily checked (see \cite{RV2})  that 
\begin{equation}\label{connection-bullets}
  \mu\bullet_\lambda\nu =\alpha_\lambda \Bigl(\Bigl( \frac{1}{\alpha_\lambda}\mu\Bigr) \bullet\Bigl(\frac{1}{\alpha_\lambda} \nu\Bigr)\Bigr) \quad\quad (\mu,\nu\in  M_e(V^K)).
\end{equation}

We now use these connections in order to study projections of  random walks on $V$ for which the involved probability
measures are contained in  $M_{e,K}^\lambda(V)$. We shall do so in discrete and continuous time and assume that either  $T:=\mathbb N_0$ or $T:=[0,\infty[$.

\begin{definition} We say that a stochastic process $(X_t)_{t\in T}$ on $V$ is a $(\lambda,K)$-compatible random walk if the starting distribution $\nu:=P_{X_0}$
  satisfies $\nu\in M_{e,K}^\lambda(V)$, and if there are (necessarily unique) probability measures $\mu_t\in  M_{e,K}^\lambda(V)$ for $t\in T$
  such that the finite-dimensional distributions of $(X_t)_{t\in T}$ satisfy
  \begin {align}\label{transition-prob-general}
    P(&X_{t_0}\in A_0,\ldots, X_{t_n}\in A_N)\\
    &= \int_V \ldots \int_V {\bf 1}_{A_0\times \ldots \times A_n}(x_0,\ldots , x_n) \> d(\delta_{x_{n-1}}*\mu_{t_n-t_{n-1}})(x_n)\ldots
     \> d(\delta_{x_{0}}* \mu_{t_1-t_{0}})(x_1)\> d\nu(x_0)
 \notag \end{align} 
  for $n\in\mathbb N$, $0=t_0<t_1<\ldots< t_n$, and Borel sets $A_0,\ldots, A_n\subset V$.
\end{definition}

By classical probability (see e.g.~\cite{Ba}), this definition means the following:\\
If $T=\mathbb N_0$ and $\nu=\delta_0$, then  $(X_t)_{t\in T}$ is just a sum of
i.i.d.~random variables on $V$ with common distribution $\mu_1\in  M_{e,K}^\lambda(V)$.
Moreover, for $T=[0,\infty[$, the process  $(X_t)_{t\in T}$ is a L\'{e}vy process with starting distribution  $\nu\in M_{e,K}^\lambda(V)$
    and associated convolution semigroup $(\mu_t)_{t\in [0,\infty[}\subset M_{e,K}^\lambda(V)$.
        The following observation is  crucial:

        \begin{lemma}\label{projected-markov}
          Let  $(X_t)_{t\in T}$ be a  $(\lambda,K)$-compatible random walk on $V$, and $p:V\to V^K$ the projection as above.
          Then $(p(X_t))_{t\in T}$ is a time-homogeneous Markov process on $V^K$ with the  finite-dimensional distributions 
 \begin {align}\label{transition-prob-projected}
    P(&p(X_{t_0})\in B_0,\ldots, p(X_{t_n})\in B_n)\\
    &= \int_{V^K} \ldots \int_{V^K} {\bf 1}_{B_0\times \ldots \times B_n}(y_0,\ldots , y_n) \> d(\delta_{y_{n-1}}\bullet_\lambda p( \mu_{t_n-t_{n-1}}))(y_n)\notag\\
    &\quad\quad\quad\quad\quad\quad\ldots
     \> d(\delta_{y_{0}}\bullet_\lambda p(\mu_{t_1-t_{0}}))(y_1)\> dp(\nu)(y_0)
   \notag  \end{align} 
  for $n\in\mathbb N$, $0=t_0<t_1<\ldots< t_n$, and Borel sets $B_0,\ldots, B_n\subset V^K$.
\end{lemma}

\begin{proof}
  We first notice that for each probability measure $\mu\in M_e(V^K)$ the mapping
  $$V^K\times \mathcal B(V^K)\to [0,1], \quad (y,B)\mapsto (\delta_y\bullet_\lambda \mu)(B)$$
  is a Markov kernel. This follows either by the very construction of $\bullet_\lambda$ above in terms of the convolution $*$, for which a corresponding result holds,
  or it can be derived from the fact that  $\bullet_\lambda$ defines a commutative hypergroup structure on $V^k$ (see \cite{RV2}),
  for which such a result holds; see the monograph \cite{BH}.

  This fact about Markov kernels ensures that the integrals of the right-hand side of (\ref{transition-prob-projected}) exist,
  and that $(p(X_t))_{t\in T}$ is a time-homogeneous Markov process on $V^K$.

  In order to check (\ref{transition-prob-projected}), we notice that (\ref{transition-prob-general}) yields that we have for the $K$-invariant sets
  $A_i:=p^{-1}(B_i)\subset V$ ($i=0,\ldots,n$) that
  \begin {align}\label{transition-prob-projected-mod}
    P(&p(X_{t_0})\in B_0,\ldots, p(X_{t_n})\in B_n)\\
    &= \int_{V} \ldots \int_{V} {\bf 1}_{A_0\times \ldots \times A_n}(x_0,\ldots , x_n) \> d(\delta_{x_{n-1}}*  \mu_{t_n-t_{n-1}})(x_n)\ldots
     \> d(\delta_{x_{0}}* \mu_{t_1-t_{0}})(x_1)\> d\nu(x_0).
  \end{align} 
  The invariance of the $A_k$ under $K$, our assumptions on the $\mu_t$, and induction on $k$ now show that in the integrations
  $ d(\delta_{x_{k-1}}*  \mu_{t_k-t_{k-1}})(x_k)$, the $x_{k-1}$ are integrated w.r.t.~some probability measures in $M_{e,K}^\lambda(V)$, and that
  step-by-step by the definition of $\bullet_\lambda$, the terms
  $\int_{V} {\bf 1}_{A_k}(x_k) \>  d(\delta_{x_{k-1}}*  \mu_{t_k-t_{k-1}})(x_k)$
can be replaced by the terms
$$\int_{V^K} {\bf 1}_{B_k}(y_k)d(\delta_{y_{k-1}}\bullet_\lambda p( \mu_{t_k-t_{k-1}}))(y_k),$$
which then leads to (\ref{transition-prob-projected}).
\end{proof}

We next compare the generators  of such Markov processes on $V^K$ for $T=[0,\infty[$ in the following setting:
 Consider some weakly continuous convolution semigroup 
$(\mu_t)_{t\in [0,\infty[}\subset M_{e,K}(V)$  of $K$-invariant probability measures on
     $V$ and  some starting distribution $\nu\in M_{e,K}(V)$. Fix some $\lambda\in V$. Then, 
 by (\ref{convo-equation}), $(e_\lambda \mu_t)_{t\in [0,\infty[}\subset M_{e,K}^\lambda(V)$ is also a semigroup.
Assume 
    from now on that this semigroup is also  weakly continuous.
Then obviously
        \begin{equation}\label{def-mu-tilde}
          \left(\tilde\mu_t:= \frac{1}{\mu_t(e_\lambda)} e_\lambda \mu_t\right)_{t\in [0,\infty[}\in M_{e,K}^\lambda(V)
          \end{equation}
        is  a weakly continuous convolution semigroup of probability measures on
        $V$. We also define the  deformed  starting distribution
        $\tilde\nu:= \frac{1}{\nu(e_\lambda)}e_\lambda\nu\in M_{e,K}^\lambda(V)$.
Now consider L\'{e}vy processes $(X_t)_{t\in [0,\infty[}$ and
              $(\tilde X_t)_{t\in [0,\infty[}$ on $V$  associated with
                  $(\mu_t)_{t\in [0,\infty[}, \nu$ and  $(\tilde \mu_t)_{t\in [0,\infty[}, \tilde\nu$
                     respectively.     
These processes are $(0,K)$- and $(\lambda,K)$-compatible respectively, and their projections 
 $(p(X_t))_{t\in [0,\infty[}$ and $(p(\tilde X_t))_{t\in [0,\infty[}$ are time-homogeneous Markov processes on $V^K$ with starting distributions
$p(\nu)$, $p(\tilde \nu)$ and transition probabilities
\begin{equation}\label{connection-transitions1}
P(p(X_t)\in A|\> p(X_s)=y)= (\delta_{y}\bullet p( \mu_{t-s}))(A),  \quad 
P(p(\tilde X_t)\in A|\> p(\tilde X_s)=y)= (\delta_{y}\bullet_{\lambda} p( \tilde\mu_{t-s}))(A)
\end{equation}
respectively for $0\le s\le t$, $y\in V^K$, and Borel sets $A\subset V^K$ where, by  Lemma \ref{modified-projection} and 
(\ref{connection-bullets}),
\begin{equation}\label{connection-transitions2}
\delta_{y}\bullet_\lambda p( \tilde\mu_{t-s})= 
\frac{1}{\alpha_\lambda(y)} \alpha_\lambda \Bigl(\delta_{y}\bullet\Bigl( \frac{1}{\alpha_\lambda} p( \tilde\mu_{t-s})\Bigr)\Bigr)=
\frac{1}{\alpha_\lambda(y)\mu_{t-s}(e_\lambda)} \alpha_\lambda (\delta_{y}\bullet p(\mu_{t-s})).
\end{equation}
This implies that the associated transition operators
\begin{equation}\label{transition-ops}
    T_t(f)(y):=\int _{V^K} f(z)\> d(\delta_{y}\bullet p( \mu_{t}))(z), \quad\quad 
    T_t^{\lambda}(f)(y)=\int _{V^K} f(z)\> d(\delta_{y}\bullet_\lambda p( \tilde\mu_{t}))(z) \quad(f\in C_0(V^K))
\end{equation}
satisfy
\begin{equation}\label{transition-ops-connection}
 T_t^{\lambda}(f)(y)=\frac{1}{\alpha_\lambda(y)\mu_{t}(e_\lambda)} T_t(\alpha_\lambda f)(y)
\end{equation}
for $t\ge 0$, $y\in V^K$, and $f\in C_b(V^K)$, where  $T_t(\alpha_\lambda f)\in C_b(V^K)$ holds   by
our assumption that for all $\mu_t$ all exponential moments exist.
Moreover, as $t\mapsto \mu_t(e_\lambda)$ is multiplicative and continuous by our assumption,
we have $\mu_t(e_\lambda)= e^{t\cdot c_\lambda}$ for all $t\ge0$ and some constant $c_\lambda\in\mathbb R$.
We thus conclude for the generators
$$Lf:= \lim_{t\to0}\frac{1}{t}(T_t(f)-f), \quad\quad L^{\lambda}f:= \lim_{t\to0}\frac{1}{t}(T_t^{\lambda}(f)-f)$$
of our transition semigroups that for
$f\in C_0(V^K)$ in the domain $D(L^\lambda)$ of $L^\lambda$,
\begin{align}\label{trafo-generators}
  L^{\lambda}f(y)&= \lim_{t\to0}\frac{1}{t}\Biggl(\frac{1}{\alpha_\lambda(y)\mu_{t}(e_\lambda)} T_t(\alpha_\lambda f)(y) -f(y)\Biggr)\notag\\
  &= \lim_{t\to0}\frac{1-\mu_{t}(e_\lambda)}{t\alpha_\lambda(y)\cdot\mu_{t}(e_\lambda)}T_t(\alpha_\lambda f)(y) +
  \frac{1}{\alpha_\lambda(y)}\lim_{t\to0}\frac{1}{t}\bigl(  T_t(\alpha_\lambda f)(y)- \alpha_\lambda(y) f(y)\bigr)\notag\\
 &= -c_\lambda f(y) +\frac{1}{\alpha_\lambda(y)} L(\alpha_\lambda f)(y).
\end{align}
This is interesting in particular for weakly continuous convolution semigroups
$(\mu_t)_{t\in [0,\infty[}\subset M_{e,K}(V)$
consisting of centered normal distributions.
In this case the 
$\tilde{\mu}_t:=\frac{1}{\mu_t(e_\lambda)} e_\lambda \mu_t \in M_{e,K}^\lambda(V)$ then are normal distributions with drifts, and the convolution semigroup $(\tilde{\mu}_t)_{t\geq0}$ clearly is weakly continuous.

 We now consider the following concrete setting which concerns  complex Cartan motion groups  where simple explicit
 formulas for the spherical functions $\alpha_\lambda$  exist. For
 the general background we refer to the monograph \cite{He} of Helgason as well to \cite{BO, Har, RV2}.

\subsection{The spherical functions of Cartan motion groups in the complex case}\label{spherical-complex-euclidean}
Let $G$ be a complex, noncompact connected 
semisimple  Lie group with finite center
with some maximal compact subgroup $K$. 
Consider the corresponding  Cartan 
decomposition $\mathfrak g= \mathfrak k\oplus \mathfrak p$ of the Lie algebra of
$G$, and consider the associated  maximal abelian subalgebra $\mathfrak a\subseteq \mathfrak p$.
The group $K$ acts on $\mathfrak p$ via the adjoint representation as
a group of orthogonal transformations w.r.t.~the Killing form as
scalar product. Let  $W$ be the Weyl group of $K$, which acts on 
  $\mathfrak a$ as a finite reflection group with  root system 
$R\subset \mathfrak a$. Then $\mathfrak a$ will be 
identified with its dual $\mathfrak a^*$ via the Killing form  $\langle\,.\,,\,.\,\rangle$. We fix some Weyl
  chamber $\mathfrak a_+$ in $\mathfrak a$ and denote the associated system of
  positive roots by  $R^+$. The closed chamber
$C:=\overline{\mathfrak a_+}$ then forms a fundamental domain for the action 
  of $W$ on $\mathfrak a$.
  
We now identify $C$ with the space $\mathfrak p^K$ of orbits where each
$K$-orbit in $ \mathfrak p$ corresponds to its unique representative 
in $C\subset\mathfrak p $.

By  Proposition IV.4.8 of
\cite{He}, the 
 $K$-invariant spherical functions on $\mathfrak p$ are given by
\begin{equation}\label{har}
\psi_\lambda(x) =  \int_K e^{i\langle\lambda,\,k.x\rangle} dk 
\quad(x\in \mathfrak p)
\end{equation}
with $\lambda$ in the complexification $\mathfrak a_\mathbb C$
 of $\mathfrak a$. 
Moreover, $\psi_\lambda\equiv \psi_\mu$ iff 
$\lambda$ and $\mu$ are in the same $W$-orbit. This is a special case of 
Harish-Chandra's integral formula for the spherical functions of a 
Cartan motion group. We now define the weighted half sum of positive roots
\begin{equation}\label{def-rho-general}\rho:= \sum_{\alpha\in R^+}\alpha\in \mathfrak a_+
\end{equation}
as in \cite{BO, RV2} (which differs by a factor $1/2$ from  the notation in \cite{He}).
Then, by Theorem II.5.35 and Cor. II.5.36 
of \cite{He} and our normalization of $\rho$,
the spherical functions $\psi_\lambda$ can be written as
\begin{equation}\label{euclidean-psi}
\psi_\lambda(x) = \frac{\pi(\rho)}{2^{|R_+|}\pi(x)\pi(i\lambda)} 
\sum_{w\in W} (\det w) e^{i \langle \lambda,\,w.x\rangle}
\end{equation}
with the fundamental alternating polynomial
$$\pi(\lambda)=\prod_{\alpha\in R^+}\langle\alpha, \lambda\rangle,$$
where $\det w$ is the determinant of the orthogonal transform $w$.
In particular,  by
Weyl's formula (see  Proposition I.5.15 of \cite{He}), 
\begin{equation}\label{weyl}\psi_{-i\rho}(x) =\, \prod_{\alpha\in R^+} 
\frac{\sinh
  \langle\alpha,x\rangle}{\langle\alpha,x\rangle}.
\end{equation}
We point out that (\ref{euclidean-psi}) and (\ref{weyl}) have to be understood in the singular cases via analytic extensions as usual.
We shall not mention this hereafter.

We now use (\ref{euclidean-psi}) and (\ref{weyl}) in order to study the projection of Brownian motions 
on 
$(\mathfrak p,\langle\,.\,,\,.\,\rangle)$ to the Weyl chamber $C$ with the methods of Subsection 2.1.
For this let $(B_t)_{t\ge0}$ be a Brownian motion on $(\mathfrak p,\langle\,.\,,\,.\,\rangle)$, i.e., w.r.t.~any orthonormal basis of
$(\mathfrak p,\langle\,.\,,\,.\,\rangle)$, the distributions $\mu_t\in M_{e,K}(\mathfrak p)$ of $B_t$ have the Lebesgue densities
$$\frac{1}{(2\pi t)^{dim \> \mathfrak p/2}}e^{-\langle x, x\rangle/(2t)} \quad\quad(t>0).$$
Now fix some $\lambda\in \mathfrak p$, and consider the function $e_\lambda(x):=e^{\langle \lambda, x\rangle}$ and the normal distributions
$$\tilde \mu_t:=\frac{1}{\mu_t(e_\lambda)} e_\lambda \mu_t\in  M_{e,K}^\lambda(\mathfrak p)$$
with the densities
$$\frac{1}{(2\pi t)^{dim \> \mathfrak p/2}}e^{-\langle x-\lambda t, x-\lambda t\rangle/(2t)} \quad\quad(t>0),$$
which belong to the Brownian motion $(B_t+\lambda t)_{t\ge0}$ on $\mathfrak p$ with drift $\lambda$.

 The Brownian motion $(B_t)_{t\ge0}$ has the  generator $\Delta/2$, and it 
 is well known  that the projection $(p(B_t))_{t\ge0}$ on $C$ then has the generator
 \begin{equation}\label{generator-projected}
   Lf(y):= \frac{1}{2}\Delta f(y)+\sum_{\alpha\in R^+}\frac{\langle \nabla f(y),  \alpha\rangle }{\langle y,  \alpha\rangle }
   \end{equation}
 for $W$-invariant functions $f\in C^{(2)}(\mathfrak a)$. This generator satisfies
 \begin{equation}\label{gen-product}
     L(f\cdot g) = f\cdot L(g)+ g \cdot  L(f)+\langle \nabla f, \nabla g\rangle
\end{equation}
 for $W$-invariant functions $f,g\in C^{(2)}(\mathfrak a)$. This follows either by a direct check, or it may be regarded as a  special case of the polarization
 of some identity for $\Delta_k(f^2)$ for general Dunkl operators $\Delta_k$ for arbitrary root systems and multiplicities $k\ge0$ in Lemma 3.1 of
 \cite{Ve}.

 We next consider the projection $(p(B_t+\lambda t))_{t\ge0}$ of the Brownian motion with drift.
 We conclude from (\ref{trafo-generators})   and  (\ref{gen-product}) with the  $W$-invariant function $g=\alpha_\lambda$ from Lemma
 \ref{modified-projection} that the generator $L^\lambda$ of this process satisfies
 \begin{align}\label{trafo-generators-special}
   L^{\lambda}f(y)&= -c_\lambda f(y) +\frac{1}{\alpha_\lambda(y)} L(\alpha_\lambda f)(y) \notag\\
   &= \Bigl( -c_\lambda+ \frac{L\alpha_\lambda(y)}{\alpha_\lambda(y)}\Bigr) f(y) +  Lf(y) +
   \frac{\langle \nabla \alpha_\lambda(y), \nabla  f(y)\rangle}{\alpha_\lambda(y)}\notag\\
    &= Lf(y) +
   \frac{\langle \nabla \alpha_\lambda(y), \nabla  f(y)\rangle}{\alpha_\lambda(y)}.
\end{align}
 Please notice here that the constant $c_\lambda$ from Subsection 2.1 satisfies $e^{tc_\lambda}=\mu_t(e_\lambda)=e^{t\|\lambda\|_2^2/2}$, and that
   $\frac{1}{2}\Delta(e_\lambda)=\|\lambda\|_2^2/2 \cdot e_\lambda$ and thus $L\alpha_\lambda =\|\lambda\|_2^2/2 \cdot \alpha_\lambda$.
   This shows that the last equation in (\ref{trafo-generators-special}) holds.
   We next conclude from Eq.~(\ref{euclidean-psi}) for  $ \alpha_\lambda=\psi_{-i\lambda}$ that
 \begin{equation}\label{log-gradient}
   \frac{1}{\alpha_{\lambda}(y)}\nabla\alpha_{\lambda}(y)= -\sum_{\alpha\in R_+}\frac{1}{\langle \alpha ,y\rangle}\alpha +
   \frac{1}{\sum_{w\in W}(\det w)  e^{\langle \lambda,w.y\rangle}}\sum_{w\in W}(\det w) e^{\langle w.\lambda,y\rangle} w.\lambda.
\end{equation}
Thus, by (\ref{trafo-generators-special}) and (\ref{generator-projected}),
\begin{equation}\label{general-generator-projected}
  L^{\lambda}f(y)=\frac{1}{2}\Delta f(y)+ \frac{1}{\sum_{w\in W}(\det w)  e^{\langle \lambda,w.y\rangle}}
  \sum_{w\in W}(\det w) e^{\langle w.\lambda,y\rangle} \langle w.\lambda,\nabla f(y) \rangle.
\end{equation}
Furthermore,  for $\lambda=\rho$, Weyl's formula (\ref{weyl}) leads to the simpler form
\begin{equation}\label{general-generator-projected-rho}
  L^{\rho}f(y)=\frac{1}{2}\Delta f(y) +\sum_{\alpha\in R^+}
  \coth(\langle \alpha, x\rangle) \langle  \alpha, \nabla f(y)\rangle.
\end{equation}
 This generator  (\ref{general-generator-projected-rho}) also appears in Section II.3 of \cite{He}.

In the remainder of this section and in Section 4  we apply the preceding results to
 examples. We begin with   the root system  $A_{N-1}$ which leads to Propositions \ref{density-a-general-drift}
and \ref{density-a-special-drift}.

\subsection{Projections of Brownian motions with drift on Hermitian matrices}\label{example-a}
The group  $K=SU(N)$ is a maximal compact subgroup of 
 the  complex, noncompact connected 
semisimple Lie  group $G=SL(N,\mathbb C)$. In the   Cartan decomposition  
$\mathfrak g= \mathfrak k \oplus\mathfrak p$ we obtain   $\mathfrak p$ as
 the additive group $H_N^0(\mathbb C)$ of
all Hermitian matrices in $H_N(\mathbb C)$ with trace $0$, on which
$SU(d)$ acts  by conjugation. 
Moreover, $\mathfrak a$
consists of all real diagonal matrices with trace 0 and will be identified
with
$$ \bigl\{x = (x_1,\ldots,x_N)\in\mathbb R^N:\> \sum_i x_i=0\bigr\}$$
on which $W$ acts as the symmetric group $S_N$ by permutations of coordinates.
We now choose  the Weyl chamber
$$C_N^{A,0}:=\{x = (x_1,\ldots,x_N)\in\mathbb R^N:\> x_1\le x_2\le\ldots\le x_N,\> \sum_i x_i=0\},$$
which parametrizes the possible spectra of matrices in $H_N^0(\mathbb C)$.
We then have the associated system of positive roots  $R^+ = \{e_j-e_i: 1\leq i<j\leq N\}$ with the standard basis
 $e_1, \ldots, e_N\in\mathbb R^N$. Furthermore, we here have the weighted half sum of positive roots
\begin{equation}\label{def-rho-a}\rho=\sum_{\alpha\in R^+}\alpha \,= \left(-N+1,-N+3,
-N+5, \ldots,N-3,N-1\right)
\in C_N^{A,0} .
\end{equation}
 The space
 $H_N(\mathbb C)$ carries the scalar product $\langle A,B\rangle:=tr(AB)$ and the associated norm $ \|.\|$.
 As
 $$H_N(\mathbb C)=H_N^0(\mathbb C)  \oplus \mathbb R \cdot I_N, \quad \text{ and } \quad
 C_N^A = C_N^{A,0} \oplus \mathbb R \cdot  (1, \ldots,1 )  \subset \mathbb  R^N,$$
 the projection $p$ from $H_N^0(\mathbb C)$ onto $C_N^{A,0}$ can be extended in the natural way to the mapping
 $p:H_N(\mathbb C) \to C_N^A$ which assigns to each Hermitian matrix its ordered spectrum.

 We now consider a Brownian motion $(B_t^H )_{ t\ge0} $ on $H_N(\mathbb C) $ as described in the introduction such that
 in particular the diagonal entries are i.i.d.~one-dimensional Brownian motions. 
We now fix some $\lambda=(\lambda_1,\ldots,\lambda_N)\in\mathbb R^N$ with $\sum_{j=1}^N \lambda_j=0$ and study the
Brownian motion $(B_t^H +t\cdot \operatorname{diag}(\lambda))_ { t\ge0} $ on $H_N(\mathbb C) $ with drift $\operatorname{diag}(\lambda)$.
If we identify vectors $x\in\mathbb R^N$ with the matrices $\operatorname{diag}(x)$, the diagonal parts of $B_t^H$ and
$B_t^H +t\cdot \operatorname{diag}(\lambda)$
have the distributions
$$d\mu_t(x)=\frac{1}{(2\pi t)^{N/2}} e^{-\|x\|^2/(2t)} dx$$
and
\begin{equation}\label{modi-a}
  d\mu_{t,\lambda}(x):= \frac{1}{(2\pi t)^{N/2}} e^{-\|x-t\lambda\|^2/(2t)} dx= e^{-\|\lambda\|^2\cdot t/2}\cdot e_{\lambda}(x) d\mu_t(x)
    \end{equation}
respectively  for $t>0$ with the functions $e_\lambda:=e_{\operatorname{diag}(\lambda)}$ from Section \ref{sub-21}. Lemma \ref{modified-projection} and
(\ref{euclidean-psi})  thus lead to the following result:

\begin{proposition}\label{density-a-general-drift-lambda}
Let  $t>0$ and $\lambda\in\mathbb R^N$ with $\sum_{j=1}^N \lambda_j=0$. Then
the ordered eigenvalues of  $(B_t^H+ t\cdot \operatorname{diag}(\lambda))_{t\ge0}$ in $C_N^A$ have the density
\begin{equation}\label{density-A-3}
  f_{N,\lambda, t}(x):=\frac{ 2^{N(N-1)/2} e^{-\|\lambda\|^2\cdot t/2}}{ (2\pi)^{N/2}  t^{N^2/2}}
  e^{-\|x\|^2/(2t)}\cdot \prod_{i<j}\frac{x_j-x_i}{\lambda_j-\lambda_i} \cdot
\sum_{w\in S_N} (\det w) e^{ \langle \lambda,\,w.x\rangle}.
  \end{equation}
\end{proposition}

\begin{proof}
This follows from Lemma \ref{modified-projection}, the  densities of the projections of the distributions  $p(\mu_{B_t^H})$ of $B_t^H$
in (\ref{density-A-1}) as well as  (\ref{euclidean-psi}),  (\ref{modi-a}), and
$$|R_+|={N(N-1)/2}   \quad\quad\text{and}\quad\quad \pi(\rho)=2^{N(N-1)/2}\cdot 2!3!\cdots (N-1)!.$$
  \end{proof}

If we apply this approach to vectors of the form $\lambda:=c\rho$ with $c>0$ and use that
 $\psi_{c\lambda}(x)=\psi_{\lambda}(cx)$
by (\ref{har}), we can apply Weyl's formula  (\ref{weyl})  instead of (\ref{euclidean-psi}). This leads to 
Proposition \ref{density-a-general-drift}.

Furthermore, if we compare the SDE (\ref{SDE1}) for $d=dim_{\mathbb R}=2$ with the generator $L^\rho$ in (\ref{general-generator-projected-rho}) for the root system
$A_{N-1}$, we see that
solutions $(X_t)_{t\ge 0}$ of (\ref{SDE1})  with d = 2 and start in 0 are in distribution equal to the processes
consisting of the ordered eigenvalues of  $(B_t^H+ t\cdot \operatorname{diag}(\rho))_{t\ge0}$. This proves
Proposition \ref{density-a-special-drift} and the following result:

\begin{corollary}\label{cor_eigval_proc}
  The  process $(\lambda_t=(\lambda_{t,1},\dots,\lambda_{t,N}))_{t\geq0}$, $\lambda_{t,1}\leq\dots\leq\lambda_{t,N}$, of the ordered eigenvalues
  of
  the Hermitian Brownian motion with drift $(B_t^H+ t\cdot \operatorname{diag}(\rho))_{t\ge0}$ is a
  solution of the SDE \eqref{SDE2} for $k=1$ with starting condition $\lambda_0=0$.
\end{corollary}

\section{Free multiplicative Brownian motions and the proof of Theorem \ref{main-theorem1}}

In this section we prove Theorem \ref{main-theorem1}.  The proof will be divided into two parts.
In the first part we  show in Theorem \ref{main-theorem1-part1} that 
the weak limits of the  normalized empirical measures $\mu_{N,t}$ in the setting of Theorem \ref{main-theorem1}
exist a.s.~for all $k\in [1/2,\infty]$ and $t>0$ and are independent from $k$.
We shall do this in the  more general setting that the  limits of the empirical starting distributions at time $0$
tend weakly to some quite arbitrary probability measure on  $\mathbb R$  while in Theorem 
\ref{main-theorem1} the empirical measures are just equal to $\delta_0$.
In the second part of the proof of Theorem \ref{main-theorem1} we then use Section 2
and describe the limit for $k=1$ via classical additive free convolutions in Theorem \ref{main-theorem1-part2}.

Let us now turn to the first result which we state in terms of free multiplicative convolutions $\boxtimes$, the spectral distribution $\mu_t$ of the free positive multiplicative Brownian motion and the exponential map $\exp_2(x):=e^{2x}$ which maps $\mathbb R$ onto $]0,\infty[$.
More precisely we will consider the free multiplicative Brownian motion $(g_t)_{t\geq0}$ with start in the identity as in the introduction, i.e. the solution to the free SDE $dg_t=i\cdot dc_t\cdot g_t$, where we integrate w.r.t. a free circular Brownian motion $(c_t)_{t\geq0}$ with variance $t$ at time $t$; see \cite{BS} for an introduction to free stochastic calculus.
In some sense $(c_t)_{t\geq0}$ can be regarded as the large $N$ operator limit of Brownian motions on $\mathbb{C}^{N\times N}$, see \cite{Ke1} for this point of view.
Our interest lies in the spectral distribution $\mu_t$ of the positive operator $g_tg_t^*$ (which agrees with the spectral distribution of the positive operator $g_t^*g_t$).
This measure was first described in \cite{B1}, and further studied in section 4.2.5 of the seminal paper \cite{B2}, see also \cite{Z,HUW} for further properties.\\
In Theorem \ref{main-theorem1-part1} below, we consider quite general starting configurations $\nu\in M^1(]0,\infty[)$ which satisfy the moment condition
\begin{equation}\label{moment_cond}
  \text{There exists some } \gamma >0\text{ with }
  s_{l}:=\int_{\mathbb{R}}x^l\,\nu(dx)\leq(\gamma l)^l\,,\quad l\in\mathbb{N}.
\end{equation}
It can be easily checked that \eqref{moment_cond} implies the  Carleman condition
\begin{equation}\label{Carleman}
		\sum_{l=1}^{\infty}s_{2l}^{-\frac{1}{2l}}=\infty\,,
	\end{equation}
which ensures that $\nu$ is determined uniquely by its moments (see e.g. \cite{A} p. 85).\\
We choose arbitrary
$x_N=(x_{N,1},\dots,x_{N,N})\in C_N^A$ for  $N\in\mathbb{N}$ such that the moments of the ``exponential'' empirical measures
satisfy
$$\frac{1}{N}\sum_{i=1}^Ne^{2lx_{N,i}}\xrightarrow{N\to\infty}s_{l} \quad\text{for}\quad l\in\mathbb{N}$$
which in particular implies that the ``exponential'' empirical measures
$\frac{1}{N}\sum_{i=1}^N\delta_{e^{2x_{N,i}}}$ tend weakly to $\nu$.
We now consider the $N$-dimensional renormalized radial Heckman-Opdam processes $\left(\tilde{X}_{t}^{(k),N}\right)_{t\geq0}$ of type A for $k\in[1/2,\infty]$
as in the introduction.
For ease of reading we drop the superscript $(k)$ in our notation for the rest of this section, i.e. appearing processes depend implicitly on the parameter $k$.
Recall that $(\tilde{X}_t^N)_{t\geq0}$ satisfies $\tilde{X}_0^{N}=x_N\in C_N^A$ and
$$
	d\tilde{X}_{t,i}^{N}
	=\begin{cases}
		\frac{1}{\sqrt{k}}\,dB_{t,i}+\sum_{j\colon j\neq i}\coth\left(\tilde{X}_{t,i}^{N}-\tilde{X}_{t,j}^{N}\right)\,dt\,,\quad&1/2\leq k<\infty\,,\\
		\sum_{j\colon j\neq i}\coth\left(\tilde{X}_{t,i}^{N}-\tilde{X}_{t,j}^{N}\right)\,dt\,,\quad&k=\infty\,,
	\end{cases}
$$
        $i\in\{1,\dots,N\}$, $N\in\mathbb{N}$.
In this situation we then have the following weak convergence result for the empirical measures of $\tilde{X}_t^{N}$, which is known for the cases $k\in\{1/2,1\}$ by \cite{Ke1} and \cite{MP}, where in the former reference a stronger type of convergence and in the latter reference a slightly weaker type of convergence was shown:

\begin{theorem}\label{main-theorem1-part1}
 	  For all $t>0$ and  $k\in[1/2,\infty]$,
$$\lim_{N\to\infty}\exp_2\left(\frac{1}{N}\sum_{i=1}^N\delta_{\tilde{X}_{t/(2N),i}^{N}}\right)=\nu\boxtimes\mu_t\;\text{weakly a.s.}$$
\end{theorem}

\begin{proof}
	We follow the proof of the corresponding result for noncompact Jacobi processes in \cite{AVW}.
  Set $Y_{t,i}^N:=\exp\left(2\tilde{X}_{t/(2N),i}^{N}\right)$ for $i\in\{1,\dots,N\}$, $t\ge0$.
  Then, by a short calculation,
	\begin{equation*}
		dY_{t,i}^N
		=\sqrt{\frac{2}{kN}}Y_{t,i}^N\,dB_{t,i}+\left[\left(\frac{N-1}{N}+\frac{1}{kN}\right)Y_{t,i}^N+\frac{2}{N}\sum_{j\colon j\neq i}\frac{Y_{t,i}^NY_{t,j}^N}{Y_{t,i}^N-Y_{t,j}^N}\right]dt
	\end{equation*}
for  $i\in\{1,\dots,N\}$ where the case $k=\infty$ is included in the obvious way.
 We want to show that the corresponding empirical measures
	$$
		\mu_{N,t}:=\frac{1}{N}\sum_{i=1}^N\delta_{Y_{t,i}^N}=\exp_2\left(\frac{1}{N}\sum_{i=1}^N\delta_{\tilde{X}_{t/(2N)}^{N}}\right)
	$$
	converge for $N\to\infty$ to  $\nu\boxtimes\mu_t$ a.s.
	We will do so by showing that the moments
	$$
		S_{N,l,t}:=\int_{\mathbb{R}}x^l\mu_{N,t}=\frac{1}{N}\sum_{i=1}^N\left(Y_{t,i}^N\right)^l$$
	 converge  to the corresponding moments of $\nu\boxtimes\mu_t$ a.s. for all
  $l\in\mathbb{N}_0$.

 Using the identity
 $$
 	2\sum_{i,j\colon j\neq i}y_iy_j\frac{y_i^{l-1}}{y_i-y_j}=\sum_{k=1}^{l-1}\sum_{i,j}y_i^ky_j^{l-k}-(l-1)\sum_{i}y_i^l\,,
 $$
 we obtain by It\^{o}'s formula
	\begin{equation*}
		dS_{N,l,t}=dM_{N,l,t}+l\left[\left(\frac{N-l}{N}+\frac{l}{kN}\right)S_{N,l,t}+\sum_{k=1}^{l-1}S_{N,k,t}S_{N,l-k,t}\right]dt\,,
	\end{equation*}
	where we set $dM_{N,l,t}:=\sqrt{\frac{2}{kN}}\frac{l}{N}\sum_{i=1}^NY_{t,i}^l\,dB_{t,i}$.
	Now define the stopping times
	$$\tau_n^l:=\inf\{t\geq0\colon S_{N,l,t}\geq S_{N,l,0}+n\}\,,\quad n,l\in\mathbb{N}\,.$$
	Then the stopped moment processes $(S_{N,l,t\wedge\tau_n^l})_{t\geq0}$ are bounded.
	Moreover, by It\^{o}'s isometry, we know that $(M_{N,l,t\wedge\tau_n^{l}})_{t\geq0}$ are $L^2$ martingales for all $n,l\in\mathbb{N}$.
	By Jensen's inequality we have the bound
	$$
		\sum_{k=1}^{l-1}S_{N,k,t\wedge\tau_n^l}S_{N,l-k,t\wedge\tau_n^l}
		\leq(l-1)S_{N,l,t\wedge\tau_n^l}\,.
	$$
	Thus
	$$
		E\left[S_{N,l,t\wedge\tau_n^l}\right]
		\leq l\left(\frac{N-l}{N}+\frac{l}{kN}+l-1\right)\int_0^tE\left[S_{N,l,s\wedge\tau_n^l}\right]\,ds\,.
	$$	
	Applying Gronwall's Lemma we get the following uniform bound in $n$:
	$$
		E\left[S_{N,l,t\wedge\tau_n^l}\right]\leq S_{N,l,0}\exp\left(l\left(\frac{N-l}{N}+\frac{l}{kN}+l-1\right)t\right)\,.
	$$
	Letting $n\to\infty$, we can deduce by monotone convergence that the expected value $E[S_{N,l,t}]$ exists for all $l\in\mathbb{N}$, $t\geq0$, and that it has the same upper bound as above.
	By It\^{o}'s isometry we see that $(M_{N,l,t})_{t\geq0}$ is an $L^2$ martingale since $E[[M_{N,l}]_t]=\frac{2l^2}{kN^2}\int_0^tE[S_{N,2l,s}]ds<\infty$, where $[M_{N,l}]_t$ denotes the quadratic variation of $M_{N,l}$ at time $t$.
	Markov- and Burkholder-Davis-Gundy inequalities yield that for some $c>0$ (independent of $N$) and all $\varepsilon>0$ we have
	\begin{align*}
		P\left(\sup_{t\in[0,T]}\left\lvert M_{N,l,t}\right\rvert>\varepsilon\right)
		\leq&\frac{c}{\varepsilon^2}E\left[M_{N,l,T}^2\right]
		=\frac{2c}{k}\left(\frac{l}{\varepsilon N}\right)^2\int_0^TE\left[S_{N,2l,t}\right]\,dt\\
		\leq&\frac{2c}{k}\left(\frac{l}{\varepsilon N}\right)^2S_{N,2l,0}\int_0^T\exp\left(2l\left(\frac{N-2l}{N}+\frac{2l}{kN}+2l-1\right)t\right)\,dt\,.
	\end{align*}
	By assumption $S_{N,2l,0}$ converges for $N\to\infty$ to $s_{2l}$, thus the term above is in $\mathcal{O}(N^{-2})$ for $N\to\infty$.
	By the Borel-Cantelli Lemma we can conclude that the limits  $\lim_{N\to\infty}S_{N,l,t}=:s_{l,t}$, $l\in\mathbb{N}_0,\,t\geq0$, exist a.s., where the convergence is locally uniformly in $t$.
	Setting $c(N,l,k):=l\left(\frac{N-l}{N}+\frac{l}{kN}\right)$, we can write the empirical moments more explicitly as follows:
	$$
		S_{N,l,t}
		=e^{c(N,l,k)t}\left[S_{N,l,0}+\int_0^te^{-c(N,l,k)s}\left(c(N,l,k)M_{N,l,s}+\sum_{k=1}^{l-1}S_{N,k,s}S_{N,l-k,s}\right)ds\right]+M_{N,l,t}\,,
	$$
	compare with (4.4)-(4.6) in \cite{AVW}.
	Letting $N\to\infty$, we thus deduce that the deterministic limit moments satisfy
	$$
		s_{l,t}
		=e^{lt}\left(s_{l}+\int_0^te^{-ls}\sum_{k=1}^{l-1}s_{k,s}s_{l-k,s}\,ds\right)\,.
	$$
	In particular $(s_{l,t})_{t\geq0}$ solves
	\begin{equation*}\label{moment_recur}
		\frac{d}{dt}s_{l,t}=l\left[s_{l,t}+\sum_{k=1}^{l-1}s_{k,t}s_{l-k,t}\right]\,,\quad s_{l,0}=s_l\,,\quad l\in\mathbb{N}\,.
	\end{equation*}
	It is well known (see the appendix in \cite{HZ} for calculations which can be used to show this fact) that this is the moment recursion of the distribution of free positive multiplicative Brownian motion, where the starting distribution $\nu$ enters by free multiplicative convolution.
	Finally we show that the moments $(s_{l,t})_{l\in\mathbb{N}_0}$ uniquely determine the measure $\nu\boxtimes\mu_t$.
	By the Carleman criterion (\cite{A} p. 85) it suffices to show that
	$$
		s_{l,t}\leq\left(e^t\gamma l\right)^l(1+t)^{l-1}\,.
	$$
	Clearly this holds for $l=1$ since $s_{1,t}=s_{1}e^t\leq e^t\gamma$.
	Moreover, we have the following estimate:
	\begin{align*}
		&\sum_{j=1}^{l-1}j^j(l-j)^{l-j}
		=2(l-1)^{l-1}+\sum_{j=2}^{l-2}j^j(l-j)^{l-j}
		\leq\left(2+4\left(\frac{l-2}{l-1}\right)^{l-1}\right)(l-1)^{l-1}\\
		\leq&(2+4/e)(l-1)^{l-1}\leq4(l-1)^{l-1}\,,
	\end{align*}
	where we used $j^j(l-j)^{l-j}\leq 4(l-2)^{l-2}$ for all $j\in\{2,\dots,l-2\}$.
	Using induction, we hence get for all $l\geq2$:
	\begin{align*}
		s_{l,t}=e^{lt}\left(s_{l}+\int_0^te^{-ls}\sum_{j=1}^{l-1}s_{j,s}s_{l-j,s}\,ds\right)
		\leq&e^{lt}\gamma^l\left(l^l+\int_0^t\sum_{j=1}^{l-1}j^j(l-j)^{l-j}(1+s)^{l-2}\,ds\right)\\
		\leq&e^{lt}(\gamma l)^l\left(1+(1+t)^{l-1}-1\right)=\left(e^t\gamma l\right)^l(1+t)^{l-1}\,.
	\end{align*}
	The moment convergence theorem yields the claim.
\end{proof}

We now turn to the setting $k=1$ and  start in $0$,
where our considerations of Section 2 give us a concrete matrix model before applying the exponential function.
This leads to another description of the limit of the empirical measures $\mu_{N,t}$ besides the one from Theorem \ref{main-theorem1-part1}.
We in turn find the following connection between the free multiplicative convolution semigroup $(\mu_t)_{t\geq0}$ and free additive convolutions
of uniform and semicircle laws:

\begin{theorem}\label{main-theorem1-part2}
	For each $N\in \mathbb N$  consider the solutions $(X_t^N:=X_t)_{t\geq0}$ of (\ref{SDE2}) for $k=1$ with start in $0\in C_N^A$.
	Then, for  $t>0$, the empirical measures $\frac{1}{N}\sum_{i=1}^N\delta_{X_{t/N,i}^{N}}$ tend weakly to $U_{t} \boxplus \mu_{sc,2\sqrt{t}}$ almost surely.
	In particular,
	$$
		\exp_2\left(U_{t} \boxplus \mu_{sc,2\sqrt{t}}\right)
		=\mu_{2t}\,,\quad t>0\,.
	$$
\end{theorem}

\begin{proof}
	Let $B_t^{H,N}:=B_t^{H}\in H_N(\mathbb{C})$ be the Hermitian Brownian motion at time $t$ as in the introduction.
	Let $\rho^N:=\rho=(-N+1,-N+3,\ldots, N-3,N-1)\in \mathbb R^N$ as in \eqref{rho-a}.
	Fix some $t>0$.
	Our first observation is that the empirical distributions
	$$
		\frac{1}{N}\sum_{i=1}^N\delta_{t(N+1-2i)/N}
	$$
	of the deterministic matrices $\frac{t}{N}\operatorname{diag}(\rho^N)$ converge for $N\to\infty$ to $U_t$.
	Also, it is well known that the empirical distributions of $B_{t/N}^{H,N}$ converge for $N\to\infty$ weakly to
        $\mu_{sc,2\sqrt{t}}$ almost surely.
	As $\sqrt{N/t}B_{t/N}^{H,N}$ has $\operatorname{GUE}(N)$ distribution, we know that asymptotically
        the matrices $B_{t/N}^{H,N}$ and $\frac{t}{N}\operatorname{diag(\rho^N)}$ are free (see Theorem 4 in section 4.2 of \cite{MS}).
	Thus the empirical distribution of $B_{t/N}^{H,N}+\frac{t}{N}\operatorname{diag}(\rho^N)$
        converges for $N\to\infty$ weakly to $U_t\boxplus\mu_{sc,2\sqrt{t}}$ almost surely.
	Moreover we know by Corollary \ref{cor_eigval_proc} that $X_{t/N}^N$ is equal
        in distribution to the ordered eigenvalues of the matrix $B_{t/N}^{H,N}+\frac{t}{N}\operatorname{diag}(\rho^N)$
        for each $N$. This implies the first claim. The second claim follows from the simple observation
        that $X_{t}^N=\tilde{X}_{t}^{N}$, $t\geq0$, since we have $k=1$, and the description of the limit of the corresponding empirical measures in Theorem \ref{main-theorem1-part1}, where our simple starting conditions $X_0^{N}=0$, $N\in\mathbb{N}$, correspond to $\nu=\delta_1$.
\end{proof}

\section{Further examples}

In this section we apply the results of Subsection 2.2 to the further infinite series of examples of complex noncompact
connected simple Lie groups with finite center which are associated with the root systems $B_N, C_N$, and $D_N$. For some data we use
Appendix C in  \cite{Kn}.

\subsection{The $B_N$-case.} For $N\ge 2$ consider the group
$G=SO(2N+1,\mathbb C)$ with  maximal compact subgroup $K=SO(2N+1,\mathbb R)$. Here $\mathfrak p$
is the vector space $Skew(2N+1,\mathbb R)$ of all skew-symmetric real matrices of dimension $2N+1$ on which 
$SO(2N+1,\mathbb R)$ acts by conjugation.
Please notice that we here suppress a possible multiplication by  $i=\sqrt{-1}$ for simplicity.

The maximal abelian subalgebra 
$\mathfrak a$  may be chosen as the vector space of all matrices  $A(x)=A(x_1,\ldots,x_N)\in Skew(2N+1,\mathbb R)$
which are formed by  $2\times 2$-blocks
of the form
$\begin{pmatrix}0& x_i\\ -x_i& 0\end{pmatrix} $
  with $x_i\in\mathbb R$ for $i=1,\ldots,N$ on the diagonal in the first $2N$ rows and columns
where all other entries (in particular in the last  row and column) are equal to $0$.  We now identify these block matrices with
$x=(x_1,\ldots, x_N)\in \mathbb R^N$ and notice that the matrices $A(x)\in \mathfrak a$ have the eigenvalues
$\pm i x_1,\ldots,\pm i x_N, 0$.

The associated Weyl group is the hyperoctahedral group
$S_N \ltimes\mathbb Z_2^N$, and the associated closed Weyl chamber
$$C_N^B:=\{0\le x_1\le\ldots\le x_N\}$$ forms a space of representatives of the $SO(2N+1,\mathbb R)$-orbits
in $Skew(2N+1,\mathbb R)$;
see e.g.~also \cite{T} for this result and the monograph \cite{M}
for related computations in the context of random matrices. The associated set of positive roots is
$$R^+=\{e_j\pm e_i:\> 1\le i<j\le N\}\cup\{e_i:\> 1\le i\le N\}$$
 with the standard
 basis $e_1,\ldots,e_N$ of $\mathbb R^N$.
 Moreover, 
 $$\rho=\left(1,3,\ldots, 2N-3, 2N-1\right),$$
 and by Weyl's formula (\ref{weyl}),
\begin{equation}\label{weyl-b}
  \psi_{-i\rho}(x) =\, \prod_{ 1\le i<j\le N} \frac{\sinh(x_j-x_i)\sinh(x_j+x_i)}{x_j^2-x_i^2} \cdot  \prod_{ 1\le i \le N}
  \frac{\sinh x_i}{x_i},
\end{equation}
Furthermore, by (\ref{generator-projected}) and (\ref{general-generator-projected-rho}), we have the generators
 \begin{equation}\label{generator-projected-b}
   Lf(y):= \frac{1}{2}\Delta f(y)+
\sum_{i=1}^N  \sum_{j:j\ne i} \Bigl(\frac{1}{y_i-y_j}+\frac{1}{y_i+y_j}\Bigr) f_{y_i}(y) +\sum_{i=1}^N \frac{1}{y_i} f_{y_i}(y)
   \end{equation}
and
\begin{equation}\label{general-generator-projected-rho-b}
  L^{\rho}f(y)=\frac{1}{2}\Delta f(y) +\sum_{i=1}^N  \sum_{j:j\ne i} \Bigl(\coth(y_i-y_j)+\coth(y_i+y_j)\Bigr) f_{y_i}(y) +\sum_{i=1}^N \coth(y_i) f_{y_i}(y),
\end{equation}
where $L^\rho$ is the projection of the generator of a Brownian motion on $G/K$.
Let $(Y_t^N)_{t\geq0}$ be the Markov process on $C_N^B$ generated by $L^{\rho}$ with start in $Y_0^N=0$.
$(Y_t^N)_{t\geq0}$ here satisfies the SDE
\begin{equation}\label{SDE4}
	dY_{t,i}^N=dB_{t,i}+\left[\coth(Y^N_{t,i})+\sum_{j\colon j\neq i}\left(\coth(Y^N_{t,i}-Y^N_{t,j})+\coth(Y^N_{t,i}+Y^N_{t,j})\right)\right]dt\,,\quad(i=1,\dots,N)\,.
\end{equation}
Results for large $N$ limits of (rescaled) empirical measures of $(Y_t^N)_{t\geq0}$ can be found in Section 5 of \cite{AVW}. We will formulate a related new result in Theorem \ref{limit_theorem_case_B} below.

We next recapitulate that the diffusions on $C_N^B$ associated with $L$ in (\ref{generator-projected-b}) with start in $0\in C_N^B$
have the densities
\begin{equation}\label{density-general-b}
  \frac{c_N^B }{t^{N^2+N/2}} e^{-\|x\|^2/(2t)} \prod_{i<j}(x_i^2-x_j^2)^{2}\cdot \prod_{i=1}^N x_i^{2}
\end{equation}
 with the normalization
\begin{equation}\label{const-b}
	c_N^B:=\frac{N!}{2^{N(N-1/2)}} \cdot\prod_{j=1}^{N}\frac{1}{j!\Gamma(\frac{1}{2}+j)}.
\notag\end{equation}
Therefore, as in the proof of Proposition \ref{density-a-general-drift} in Subsection 2.3, we obtain from Lemma
\ref{modified-projection} and (\ref{weyl-b}):

\begin{lemma}\label{density-b-general-drift} Let $(B_t)_{t\ge0}$ be a Brownian motion on  $Skew(2N+1,\mathbb R)$ with start in $0$ and
  $(X_t:=B_t+t\cdot A(\rho))_{t\ge0}$ the associated Brownian motion with drift $ A(\rho)$. Then the diffusion
  on $C_N^B$, which describes the eigenvalues of $(X_t)_{t\ge0}$ (up to multiplication by $i$ as described above), has the generator
  $L^{\rho}$ from (\ref{general-generator-projected-rho-b}) and, for $t>0$, the densities
  $$  \frac{ c_N^Be^{-\|\rho\|^2 t/2}}{t^{N^2+N/2}} e^{-\|x\|^2/(2t)} \prod_{i<j}(x_i^2-x_j^2)\sinh(x_i-x_j)\sinh(x_i+x_j)\cdot \prod_{i=1}^N x_i \cdot\sinh x_i$$
 with $\lvert\lvert\rho\rvert\rvert^2=N(2N-1)(2N+1)/3$.
 \end{lemma}

The preceding results  lead to a free limit for Brownian motions on $SO(2N+1,\mathbb C)/SO(2N+1,\mathbb R)$ for
$N\to\infty$ similar to Theorems \ref{main-theorem1-part2} and \ref{main-theorem1}:

\begin{theorem}\label{limit_theorem_case_B}
  For each $N\in \mathbb N$ consider the solution $(Y_t^N)_{t\geq0}$ of (\ref{SDE4}) with start in $0\in C_N^B$.
Then, for all $t>0$, the  normalized empirical measures
	\begin{equation}
		\mu_{N,t}:=\frac{1}{N}\sum_{i=1}^N\delta_{Y^N_{t/(2N),i}}
	\end{equation}
        tend weakly to $\lvert U_{t}\boxplus\mu_{sc,2\sqrt{t}}\rvert$ almost surely where $|\mu|$ denotes the push forward of the measure $\mu$
        under $x\mapsto|x|$.
\end{theorem}

\begin{proof}
  Let $s>0$, and let $(Z^N_t)_{t\geq0}$ be the Markov process on $C_N^B$ with generator  $L$ from \eqref{generator-projected-b}
  with start in $Z^N_0:=\frac{s}{2N}\rho$. By the form of $\rho$, the empirical measures of the starting configuration converge
  weakly to the uniform distribution on $[0,s]$, i.e.,
	$$
		\lim_{N\to\infty}\frac{1}{N}\sum_{i=1}^N\delta_{\frac{s}{2N}\rho_i}
		=\operatorname{Unif}([0,s])\,.
	$$
	Moreover the unique even measure $\mu$ on $\mathbb{R}$ with $|\mu|=\operatorname{Unif}([0,s])$ is given by $U_s=\operatorname{Unif}([-s,s])$.
	We thus conclude from Theorem 5.2 in \cite{AVW} that
	$$
		\lim_{N\to\infty}\frac{1}{N}\sum_{i=1}^N\delta_{Z_{t/(2N)}^N}
		=\lvert U_s\boxplus\mu_{sc,2\sqrt{t}}\rvert\quad\text{weakly a.s.}\,,\;t>0\,,
	$$
	where we used the fact that a Bessel process of type B at time $t$ with start in $y$ is equal in distribution to a Bessel process of type B at time $t/(2N)$ with start in $y/\sqrt{2N}$ scaled in space by the factor $\sqrt{2N}$.
	The claim now follows from Lemma \ref{density-b-general-drift}, since for the choice $s=t$ we know that $Y_{t/(2N)}^N$ and $Z_{t/(2N)}^N$ are equal in distribution for all $N\in\mathbb{N}$.
\end{proof}

\subsection{ The $D_N$-case} For $N\ge 2$ consider the group
$G=SO(2N,\mathbb C)$ with  maximal compact subgroup $K=SO(2N)$. Here $\mathfrak p=Skew(2N,\mathbb R)$, and
$\mathfrak a$ may be chosen as the vector space of all matrices $A(x)$  with  $2\times 2$-blocks
of the form
$\begin{pmatrix}0& x_i\\ -x_i& 0\end{pmatrix} $
  with $x_i\in\mathbb R$ for $i=1,\ldots,N$ on the diagonal and all other entries equal to $0$.
Again we identify  $\mathfrak a=\mathbb R^N$. Here the Weyl chamber 
$$C_N^D=\{x\in\mathbb R^N:\> x_N\ge x_{N-1}\ge\cdots\ge x_{2}\ge |x_1|\}$$
forms a space of representatives of the $K$-orbits in $Skew(2N,\mathbb R)$ with the associated
set
$$R^+=\{e_j\pm e_i:\> 1\le i<j\le N\}$$
of positive roots and
 $$\rho=\left(0,2,\ldots, 2(N-3), 2(N-1)\right),$$
 and by Weyl's formula (\ref{weyl}),
\begin{equation}\label{weyl-d}
  \psi_{-i\rho}(x) =\, \prod_{ 1\le i<j\le N} \frac{\sinh(x_j-x_i)\sinh(x_j+x_i)}{x_j^2-x_i^2}.
\end{equation}
Furthermore, by (\ref{generator-projected}) and (\ref{general-generator-projected-rho}), we have the generators
 \begin{equation}\label{generator-projected-d}
   Lf(y):= \frac{1}{2}\Delta f(y)+
\sum_{i=1}^N \sum_{j:j\ne i}\Bigl( \frac{1}{y_i-y_j}+ \frac{1}{y_i+y_j}\Bigr) f_{y_i}(y) 
   \end{equation}
and
\begin{equation}\label{general-generator-projected-rho-d}
  L^{\rho}f(y)=\frac{1}{2}\Delta f(y) +\sum_{i=1}^N  \sum_{j:j\ne i} \Bigl(\coth(y_i-y_j)+\coth(y_i+y_j)\Bigr) f_{y_i}(y).
\end{equation}
The diffusions on $C_N^D$ associated with (\ref{generator-projected-d}) with start in $0\in C_N^D$
have the densities
\begin{equation}\label{density-general-d}
  \frac{c_N^D }{t^{N^2-N/2}} e^{-\|x\|^2/(2t)} \prod_{i<j}(x_i^2-x_j^2)^{2}
\end{equation}
 with the normalization
\begin{equation}\label{const-d}
 c_N^D:=
 \frac{N!}{2^{N(N-3/2)+1}} \cdot\prod_{j=1}^{N}\frac{1}{j!\Gamma(j-\frac{1}{2})}.
\notag\end{equation}
Therefore, as in the case above:

\begin{lemma}\label{density-d-general-drift} Let $(B_t)_{t\ge0}$ be a Brownian motion on  $Skew(2N,\mathbb R)$ with start in $0$ and
  $(X_t:=B_t+t\cdot A(\rho))_{t\ge0}$ the associated Brwwnian motion with drift $ A(\rho)$. Then the diffusion
  on $C_N^D$, which describes the eigenvalues of $(X_t)_{t\ge0}$ (up to multiplication by $i$), has the generator
  $L^{\rho}$ from (\ref{general-generator-projected-rho-d}) and, for $t>0$, the densities
  $$  \frac{ c_N^De^{-\|\rho\|^2 t/2}}{t^{N^2-N/2}} e^{-\|x\|^2/(2t)} \prod_{i<j}(x_i^2-x_j^2)\sinh(x_i-x_j)\sinh(x_i+x_j)$$
  with $\lvert\lvert\rho\rvert\rvert^2=2N(N-1)(2N-1)/3$.
  \end{lemma}

Free limits in the D-case lead to the same results as for the B-case in Theorem \ref{limit_theorem_case_B}
 
\subsection{ The $C_N$-case} For $N\ge 2$ consider
$$G=Sp(N,\mathbb C):=\{A\in SL(2N,\mathbb C):\>\>   A^tJ_N A=J_N\} \quad\text{with the $2N\times 2N$-matrix}\quad
J_N=\begin{pmatrix}0& I\\ -I& 0\end{pmatrix} $$
with   maximal compact subgroup $K=Sp(N)=\{A\in GL(N,\mathbb H):\>\>   A^* A=I\} $.
In this case, we again can identfy   $\mathfrak a$ with $\mathbb R^N$ and obtain the Weyl chamber
 $C_N^B$ and the Weyl group $W$ as in the $B_N$-case.
 We then have the positive root system 
$$R^+=\{e_j\pm e_i:\> 1\le i<j\le N\}\cup\{2e_i:\> 1\le i\le N\},$$
the vector  $\rho=(2,4,\ldots, 2N-2,2N)$.
 and by Weyl's formula (\ref{weyl}),
\begin{equation}\label{weyl-c}
  \psi_{-i\rho}(x) =\, \prod_{ 1\le i<j\le N} \frac{\sinh(x_j-x_i)\sinh(x_j+x_i)}{x_j^2-x_i^2} \cdot  \prod_{ 1\le i \le N}
  \frac{\sinh (2x_i)}{2x_i}.
\end{equation}
Furthermore, by (\ref{generator-projected}) and (\ref{general-generator-projected-rho}), we have the generator $L$ from (\ref{generator-projected-b}),
which does not change here,
as well as
\begin{equation}\label{general-generator-projected-rho-c}
  L^{\rho,C}f(y)=\frac{1}{2}\Delta f(y) +\sum_{i=1}^N  \sum_{j:j\ne i} \Bigl(\coth(y_i-y_j)+\coth(y_i+y_j)\Bigr) f_{y_i}(y) +2\sum_{i=1}^N \coth(2y_i) f_{y_i}(y).
\end{equation}
If we use the densities (\ref{density-general-b}) associated with  the generators $L$ from (\ref{generator-projected-b}),
 the methods of the proof of Proposition \ref{density-a-general-drift} in Subsection 2.3 lead to:

 \begin{lemma}\label{density-c-general-drift} Let $(X_t)_{t\ge0}$ be a diffusion on $C_N^B$ with start in $0$, which belongs to the generator
 $L^{\rho,C}$. Then for $t>0$, $X_t$ has the density
  $$  \frac{ c_N^Be^{-\|\rho\|^2 t/2}}{t^{N^2+N/2}} e^{-\|x\|^2/(2t)} \prod_{i<j}(x_i^2-x_j^2)\sinh(x_i-x_j)\sinh(x_i+x_j)\cdot \prod_{i=1}^N \frac{x_i}{2} \cdot\sinh (2x_i)$$
 with $\lvert\lvert\rho\rvert\rvert^2=2N(N+1)(2N+1)/3$.
 \end{lemma}

 Again, free limits in this C-case lead to the same results as for the B-case above.

 \begin{remarks}
\begin{enumerate}
\item[\rm{(1)}] All results above for the particular drift $\rho$ can be also stated for arbitrary vectors $\lambda$ by using (\ref{euclidean-psi})
 instead of Weyl's formula (\ref{weyl}) similar to  Lemma \ref{density-a-general-drift-lambda}.
\item[\rm{(2)}] Consider the groups $G=SU(N,M)$ with maximal compact subgroups $K=SU(N)\times SU(M)$  over $\mathbb C$ for $M\ge N$, i.e.,
  the symmetric spaces  $G/K$
  are  noncompact Grassmann manifolds. Here we can identify $\mathfrak p$ with $\mathbb C^{N,M}$,
  and $\mathfrak a$ with $\mathbb R^N$, and we have the root system $BC_N$; see e.g. Section 6 of \cite{BO}.
In these cases there
exist  determinantal formulas for the spherical functions due to Berezin and  Karpelevich \cite{BK}
(see also \cite{Ho} and references there)
in terms of one-dimensional Jacobi functions; see  \cite{Ko}
for  these functions. Moreover, there are related determinantal formulas for the associated Euclidean spherical functions
where  Bessel functions instead of Jacobi functions appear in these determinants; see \cite{Me, BO}.
Now consider  the projections of Brownian motions on  $\mathbb C^{N,M}$
  with drifts to the associated Weyl chambers $C_N^B$. Then Section 2 leads to explicit
  determinantal formulas for the densities on   $C_N^B$
  similar to those in Lemma \ref{density-a-general-drift-lambda}. However, these formulas do not describe corresponding
  densities of Brownian motions on $G/K$ as it was the case for all other examples above.
\item[\rm{(3)}] Consider some  distribution  $P\in M^1(C)$
  which appears in Proposition \ref{density-a-general-drift} or  Lemmas \ref{density-b-general-drift}
  \ref{density-d-general-drift}, \ref{density-c-general-drift}, or in the preceding remark on the corresponding Weyl chamber
  $C\subset \mathbb R^ N$. Consider the associated 
   Weyl group invariant measure $\tilde P:=\frac{1}{|W|}\sum_{w\in W} w(P) \in M^1(\mathbb R)$.
  It would be interesting to determine the $k$-dimensional 
  marginal distributions ($k=1,\ldots,N-1$) and level spacing distributions
  of $\tilde P$, which are well known, when no drift appears; see the monograph \cite{M}.
\end{enumerate}
   \end{remarks}


\begin{thebibliography}{999}
  
\bibitem[A]{A} N.I. Akhiezer, The Classical Moment Problems and Some Related Questions in Analysis.
 Engl. Translation, Hafner Publishing Co., New York, 1965. 

\bibitem[AGZ]{AGZ} G.W. Anderson, A. Guionnet, O. Zeitouni, 
  An Introduction to Random Matrices. Cambridge University Press, Cambridge 2010.


\bibitem[AVW]{AVW} M. Auer, M. Voit, J.H.C. Woerner,  Wigner- and Marchenko-Pastur-type limits for Jacobi processes.
  \textit{J. Theor. Probab.}  37 (2024), 1674-1709. 

\bibitem[Ba]{Ba} H. Bauer, Probability Theory. de Gruyter, Berlin (2002).

\bibitem[BK]{BK} F. A. Berezin, F. I. Karpelevich, Zonal spherical functions and Laplace operators on some symmetric spaces,
  \textit{Dokl. Akad. Nauk SSSR} 118(1) (1958), 9-12.

\bibitem[BO]{BO} S. Ben Saïd, B. Ørsted, Analysis on flat symmetric spaces. \textit{J. Math. Pures Appl.} 84 (2005), 1393–1426.

\bibitem[B1]{B1} P. Biane,
Free Brownian motion, free stochastic calculus and random matrices. In:
Free probability theory (Waterloo, ON, 1995), 1-19, \textit{Fields Inst. Commun.}, 12, Amer. Math. Soc., Providence, RI, 1997.

\bibitem[B2]{B2} P. Biane, 
Segal-Bargmann transform, functional calculus on matrix spaces and the theory of semi-circular and circular systems. 
\textit{J. Funct. Anal.} 144 (1997),  232-286.

\bibitem[B3]{B3} P. Biane,
Processes with free increments.
\textit{Math. Z.} 227 (1998), 143–174.

\bibitem[BS]{BS} P. Biane, R. Speicher, Stochastic calculus with respect
  to free Brownian motion and analysis on Wigner space, \textit{Probab. Theory Rel. Fields} 112 (1998), 373-409.

  \bibitem[BH]{BH} W. Bloom, H. Heyer, \textit{Harmonic Analysis of 
    Probability Measures on Hypergroups.} De Gruyter-Verlag, Berlin, 1994.

\bibitem[CGY]{CGY} O. Chybiryakov, L. Gallardo, M. Yor, Dunkl processes and their
 radial parts relative to a root system. In:
 P. Graczyk et al. (eds.), Harmonic and stochastic analysis of Dunkl processes, pp. 113-198. Hermann, Paris 2008.




\bibitem[DV]{DV} J.F. van Diejen, L. Vinet, Calogero-Sutherland-Moser Models.
  CRM Series in Mathematical Physics, Springer, Berlin, 2000.


 \bibitem[F]{F}  P. Forrester, Log Gases and Random Matrices, London Mathematical Society, London, 2010.

  
   
\bibitem[GM]{GM} P.  Graczyk, J. Malecki, Strong solutions of non-colliding particle systems.
  \textit{ Electron. J. Probab.} 19 (2014), 21 pp.


\bibitem[HZ]{HZ} Hao-Wei  Huang,  Ping Zhong, 
On the supports of measures in free multiplicative convolution semigroups.
 \textit{Math. Z.} 278 (2014),  321-345.


 \bibitem[Har]{Har}  Harish-Chandra, Representations of semisimple Lie groups, III,  \textit{Trans. Amer. Math. Soc.} 76 (1954), 234–
253.

\bibitem[HUW]{HUW} T. Hasebe, Y. Ueda, J.-C. Wang,
Log-unimodality for free positive multiplicative Brownian motion.
\textit{Coll. Math.} 169 (2022), 209-226.

\bibitem[HO]{HO} G. Heckman, E. Opdam, Jacobi polynomials and hypergeometric functions associated with root systems.
  In: Encyclopedia of Special Functions, Part II: Multivariable Special Functions, eds. T.H. Koornwinder, J.V. Stokman,
  Cambridge University Press, Cambridge, 2021.

\bibitem[HS]{HS} G. Heckman, H. Schlichtkrull, 
Harmonic Analysis and Special Functions on Symmetric Spaces, Part I.  
Perspectives in Mathematics, Vol. 16, Academic Press, 1994.



\bibitem[He]{He} S. Helgason, Groups and Geometric Analysis.
Academic Press, 1984.


\bibitem[Ho]{Ho} B. Hoogenboom,
  Spherical functions and invariant differential operators on complex Grassmann manifolds,  \textit{Ark. Mat.}
  20 (1982), 69-85.
  
\bibitem[IS]{IS}  J. R. Ipsen, H. Schomerus,
  Isotropic Brownian motions over complex fields as a solvable model for May–Wigner stability analysis.
\textit{J. Physics A}  49, Number 38 (2016),  385201.



\bibitem[Ke1]{Ke1} T. Kemp,
The large-N limits of Brownian motions on $GL_N$.
\textit{Int. Math. Res. Notices} 13 (2016), 4012-4057.

\bibitem[Ke2]{Ke2}  T. Kemp, 
Heat kernel empirical laws on $U_N$ and 
and $GL_N$.
\textit{J. Theor. Probab.} 30 (2017),  397-451.

\bibitem[Kl1]{Kl1} A.  Klyachko, Stable bundles, representation theory, and
    Hermitian operators,  \textit{Selecta Math.} (new Series) {4} (1998),
419-445. 


\bibitem[Kl2]{Kl2} A.  Klyachko, Random walks on symmetric spaces
  and inequalities for matrix spectra, \textit{Linear Algebra Appl.}
  {319} (2000), 37--59.

\bibitem[Kn]{Kn} A.W.  Knapp, Lie Groups beyond an Introduction.
  Birkh\"auser, Boston, 1996.

\bibitem[Ko]{Ko} T. Koornwinder, Jacobi functions and 
analysis of noncompact semisimple
Lie groups. In: R. A. Askey et al. (eds.), Special Functions:
Group Theoretical Aspects
and Applications. Dordrecht-Boston-Lancaster: D. Reidel
Publishing Company 1984, pp. 1-85.

\bibitem[MP]{MP} J. Ma\l{}ecki, J. L. P{\'e}rez, Universality classes for general random matrix flows.
  \textit{Ann. Inst. Henri Poincar\'{e}}, Prob.  Stat. 58 (2022),  722-754.

\bibitem[Me]{Me}   C. Meaney, The inverse Abel transform for $SU(p, q)$, \textit{Ark. Mat.} 24 (1985), 131–140.

\bibitem[M]{M} M. Mehta, Random matrices (3rd ed.), Elsevier/Academic Press, Amsterdam, 2004.
  
\bibitem[MS]{MS} J. Mingo, R. Speicher: Free Probability and Random Matrices, Fields Institute Monographs, Springer, 2017
  
\bibitem[NS]{NS} A. Nica, R. Speicher, Lectures on the Combinatorics of Free Probability Theory, Cambridge University Press, Cambridge, 2006.
  
\bibitem[NRW]{NRW} J. R. Norris,  L. C. G. Rogers, D. Williams,
Brownian motions of ellipsoids. 
\textit{Trans. Am. Math. Soc.} 294 (1986), 757-765.


\bibitem[P]{P} P.E. Protter, Stochastic Integration and Differential Equations. A New Approach.
  Springer, Berlin, 2003.




\bibitem[RV1]{RV1} M. R\"osler, M. Voit, Markov processes related with Dunkl operators.
\textit{Adv. Appl. Math.}  21 (1998), 575-643.

\bibitem[RV2]{RV2} M. R\"osler, M. Voit, $SU(d)$-biinvariant random walks on $SL(d,C)$
and their Euclidean counterparts. 
 \textit{Acta Appl. Math.} 90 (2006),  179-195. 


  

\bibitem[Sch1]{Sch1} B. Schapira,
The Heckman-Opdam Markov processes.
 \textit{Probab. Theory Rel. Fields} 138 (2007), 
 495-519.

\bibitem[Sch2]{Sch2} B. Schapira, Contribution to the hypergeometric function theory of Heckman and Opdam: 
  sharp estimates, Schwarz space, heat kernel.  \textit{Geom. Funct. Anal.} 18 (2008), 222-250.

\bibitem[TK]{TK} Y. Takahashi, M. Katori, Noncolliding Brownian motion with drift and time-dependent Stieltjes-Wigert determinantal point process. \textit{Journal of Mathematical Physics} 53, No. 10, 103305 (2012).

\bibitem[T]{T} G. Thompson, Normal forms for skew-symmetric matrices and Hamiltonian systems with first integrals linear in momenta. 
  \textit{Proc. Amer. Math. Sic.} 104 (1988), 910-916.

    
 \bibitem[Ve]{Ve}  A. Velicu, 
Sobolev-type inequalities for Dunkl operators.
\textit{J. Funct. Anal.} 279, No. 7, Article ID 108695, 36 p. (2020). 
  
 \bibitem[V]{V} M. Voit,    Positive
 characters on commutative
hypergroups and some  applications,  \textit{ Math. Z.} 198 (1988), 405-421.

 \bibitem[VW1]{VW1} M. Voit, J.H.C. Woerner,  Limit theorems for Bessel and Dunkl processes of large dimensions and free convolutions.
 \textit{Stoch. Proc. Appl.} 143 (2022), 207-253.

\bibitem[VW2]{VW2} M. Voit, J.H.C. Woerner, The differential equations associated with Calogero-Moser-Sutherland
  particle models in the freezing regime. \textit{Hokkaido Math. J.} 51 (2022), 153--174.
  
\bibitem[Z]{Z} P. Zhong,
On the free convolution with a free multiplicative analogue of the normal
distribution.
\textit{J. Theoret. Probab.} 28 (2015), no. 4, 1354-1379.

\end{thebibliography}
\end{document}